\documentclass{amsart}
\date{}

\title[Curved flats]{Curved flats in symmetric spaces}

\author[D. Ferus]{Dirk Ferus}
\address{Fachbereich Mathematik, MA 8-3 \\
Technische Universit\"at Berlin \\
Strasse des 17. Juni 136\\ 10622 Berlin, Germany}
\email[Dirk Ferus]{ferus@sfb288.math.tu-berlin.de}

\author[F. Pedit]{Franz Pedit}
\address{Department of Mathematics \\
University of Massachusetts \\
Amherst, MA 01003}
\email[Franz Pedit]{franz@gang.umass.edu}
\thanks{Second author was partially supported by NSF grant DMS-9205293,
the SFB 288 at Technische Universit\"at Berlin and the
Graduierten Kolleg at Humboldt Universit\"at Berlin.}

\newcommand{\C}{{\mathbb C}}

\newcommand{\Rr}{{\mathbb R}}

\newcommand{\p}{{\mathfrak p}}
\renewcommand{\a}{{\mathfrak a}}

\newcommand{\Ad}{\operatorname{Ad}}
\renewcommand{\d}{\operatorname{d}}
\renewcommand{\dim}{\operatorname{dim}}
\newcommand{\rank}{\operatorname{rank}}

\newcommand{\R}[1]{\mbox{${{\mathbb R}^#1}$}}
\newcommand{\Lorentz}[1]{\mbox{${\mathbb R}_1^{#1}$}}
\newcommand{\Lor}{\Lorentz{2m+1}}

\newtheorem{mytheorem}{Theorem}

\newtheorem{mylemma}[mytheorem]{Lemma}

\theoremstyle{definition}
\newtheorem*{mydefinition}{Definition}
\newtheorem{example}{Example}

\theoremstyle{remark}
\newtheorem{remark}{Remark}

\numberwithin{equation}{section}

\begin{document}
\maketitle

\section{Introduction}
Over the past years many problems in classical surface theory and also
more generally, sub\-mani\-fold theory, have been linked to certain types of
completely integrable non-linear PDE (soliton equations). In certain
cases this has already been known to geometers in the last century
but the explicit construction of solutions by the finite gap
integration scheme has only been achieved recently. The crucial ingredient
is to rewrite the equations at hand as zero-curvature equations
involving a {\em spectral} parameter --- usually given by a non-trivial
geometrical
deformation of the surface. The finite gap theory then provides special
solutions to those zero-curvature equations parameterized by certain
algebraic curves, the {\em spectral curves}, which in certain cases account for
all solutions. A standard
example is that of (pluri) harmonic maps into symmetric spaces \cite{BFPP}
which
contains many examples of surface theory -- constant mean and Gauss curvature
surfaces, Willmore surfaces, minimal immersions --- as special cases.
Recently also the isometric immersion equations for maps between
space forms (with flat normal bundle) \cite{FP}, conformally flat
$3$-folds in the $4$-sphere \cite{UJ, UJ1} and isothermic surfaces \cite{BJPP}
have been shown to fit into this approach.

In this paper we introduce a natural class of maps into symmetric spaces,
{\em curved flats},
which, in a certain sense,  contains all the above examples as special cases.
A curved flat is simply a map into a symmetric space
which is tangent at each point to a flat of the symmetric space, i.e.
each tangent space is abelian. Since curved flats are intrinsically flat
(Theorem~\ref{th:flat}) they
may also be viewed as natural analogues of developable surfaces in $3$-space.

Important examples of curved flats arise from isometric immersions of space
forms
into space forms via their Gauss maps (Section~\ref{sec:examples}).
The developing isometry identifies (locally) the curved flat with a Cartan
subalgebra of the Grassmanian (in which the Gauss map takes values in).
The roots of the Cartan subalgebra are mapped via this isometry into
principal curvature coordinates of the corresponding isometric immersion
(Section~\ref{sec:geometry}). Further examples of curved flats
come from isothermic surfaces \cite{BJPP} and conformally flat
$3$-folds in the $4$-sphere \cite{UJ, UJ1}. Of course, every curve in a
symmetric space
is a curved flat.

Besides their geometric relevance curved flats are also interesting from
the integrable systems point of view:  scaling the derivative of a curved
flat gives rise to a non-trivial deformation which obviously preserves the
curved flat condition (Section~\ref{sec:spectral}). Therefore one can write
the curved flat condition as a zero-curvature equation
on a Lie algebra valued $1$-form involving the deformation -- {\em spectral} --
parameter
linearly. Introducing the necessary loop algebras in
Section~\ref{sec:integrate}
we then apply the finite type integration scheme to construct solutions to
the curved flat equations from a hierarchy of commuting algebraically
integrable ODE in Lax form. This gives a space of solutions which has the
functional
dimension
of finitely many functions of a single variable.  At least
in principle, each such solution can be expressed by theta functions
on some algebraic curve. In the rather instructive simple example of a curve in
the
$2$-sphere (which we will not persue in this paper)
the spectral parameter scales the speed of the curve and at the same time
scales the (geodesic) curvature inverse proportional. The lowest type curves
give the geodesics (spectral curve is the Riemann sphere) and elastic
curves (spectral curve is elliptic).

We would like to thank Fran Burstall, Ottmar Loos and Ulrich Pinkall for
stimulating discussions during the preparation of this paper. Finally,
the second author would like to thank the members of the SFB 288 at TU-Berlin
for their interest and hospitality.

\section{Curvature isotropic immersions}

Let $N = G/K$ be a semisimple symmetric space endowed with the
Killing metric,
$M$ a manifold of dimension $m$, and  $f: M\rightarrow N$
a smooth map. Since most of our considerations are local, we assume $M$
to be contractible unless otherwise stated.

\begin{mydefinition}\label{def:flat}
 $f$ is called {\em curvature isotropic } or a
{\em curved flat}, if we have
\[f^* R^N=0 \]
as a 2-form on $M$ with values in the bundle $End(f^* TN)$.
Here $R^N$ denotes the curvature tensor of $N$.

Moreover, if $f$ is an immersion which induces a nondegenerate (not necessarily
positive definite) metric on M,
we call $f$ a {\em regular } curved flat.

\end{mydefinition}

We give a Lie algebraic formulation of this definition, which will also
explain the term {\em curved flat}.
Let
\[\mathfrak{g}=\mathfrak{k}+\mathfrak{p}.\]
be the symmetric decomposition of $G/K$.
Let
$F: M \rightarrow G$  be a lift
of \(f\) so that $\pi\circ F=f$ where $\pi:G\to G/K$, $\pi(g)=[g]$, is the
coset projection.
We call such a lift a {\em framing {\rm of} f}. Then
\begin{equation}
F^{-1}\d F=A=A_0+A_1, \label{decomp}
\end{equation}
where $A$ is a 1-form on $M$ with values in $\mathfrak{g}$,
and $A_0$  resp. $A_1$  denote the $\mathfrak{k}$ resp. $\mathfrak{p}$
parts of $A$.
$A$ satisfies the Maurer-Cartan equation or
integrability condition
\begin{equation}
\d A+\frac{1}{2}[A\wedge A]=0,\label{MC}
\end{equation}
where
\[
[A\wedge B](X,Y)=[A(X),B(Y)]-[A(Y),B(X)]\,.
\]
The $\mathfrak{p}$-space is closed under the triple Lie bracket
and canonically identified
with the tangent space of $N$ at the origin $[e]=K$.
Then, if \(f(p)=[e], F(p)=e\), and \(X\in T_pM\), we obtain

\begin{eqnarray*}
f_* (X) &= &\pi_* F_* (X)=\pi_* (FA(X))\\
                &=&\pi_* (A_0(X)+A_1(X))=\pi_* A_1(X) \\
                &=&A_1(X),
\end{eqnarray*}

and
\begin{eqnarray*}
f^* R^N(X,Y) &=&- \frac{1}{2}[[A_1(X),A_1(Y)],...]  \\
&=&- \frac{1}{2}[[A_1\wedge A_1](X,Y),...]\,.
\end{eqnarray*}

\begin{mylemma}\label{lem:abel}
The subspace $\mathfrak{a}:=A_1(T_pM)$ is abelian:
\begin{equation}
[A_1\wedge A_1]=0. \label{flat}
\end{equation}
\end{mylemma}

\begin{proof}
The assumption implies
\begin{equation}
[[\mathfrak{a},\mathfrak{a}],\mathfrak{p}]=0. \label{abel}
\end{equation}
Note that
\[
\mathfrak{g_1}:=[\mathfrak{p},\mathfrak{p}] + \mathfrak{p} \mbox{ is an ideal
of }%
\mathfrak{g}.
\]
Since $\mathfrak{g}$ is semisimple, so is $\mathfrak{g_1}$, and by
(\ref{abel}),
$[\mathfrak{a},\mathfrak{a}]$ is in the center of $\mathfrak{g_1}$. Hence
$[\mathfrak{a},\mathfrak{a}]=0$.
\end{proof}

Conversely, a $\mathfrak{g}$-valued 1-form $A$ that satisfies
the integrability condition (\ref{MC}) together with (\ref{flat})
gives rise to a curvature isotropic map f.

The Lemma shows that a curved flat is at each point tangent to a flat
totally geodesic submanifold of $N$.

In view of (\ref{flat}) the $\mathfrak{k}$- and $\mathfrak{p}$-part of
(\ref{MC})
reduce to
\begin{eqnarray}
\d A_0+\frac{1}{2}[A_0\wedge A_0]&=&0 \label{MCk}
\\
\d A_1+[A_0\wedge A_1]&=&0. \label{MCp}
\end{eqnarray}
Note that \(A_0:TM \rightarrow \mathfrak{k}\) is the canonical connection of
the principal \(K\)-bundle \(f^{*}(G) \rightarrow M\), and this is
flat by equation (\ref{MCk}).

\section{Examples}\label{sec:examples}

\begin{example}\label{ex:curves}
Any curve in a symmetric space, in particular in the rank 1 space $S^2$ is
curvature isotropic. This rather trivial example nicely displays
several features of the general theory in a very explicit form. We shall
come back to this elsewhere.
\end{example}

\begin{example}\label{ex:isomequal}
 Let $M^m$ be of constant curvature $c>0$, and let
\[
\phi : M^m \rightarrow S:=S^{n+1}(c)\subset \R{{n+2}}
\]
be an isometric immersion into the sphere of equal curvature $c$.
The Gauss map of $\phi$
\[
f=\gamma_\phi : M\rightarrow G_{m+1}(\R{{n+2}})=
\frac{{\bf O}(n+2)}{{\bf O}(m+1)\times {\bf O}(n+1-m)}
\]
assigns to each point of M the orthogonal complement of
its normal space in
the sphere $S$ considered as linear subspace of \R{{n+2}}.  If
$\phi$ has flat normal bundle $V^{\perp}$, we claim that
\[
\gamma_\phi \; \mbox{is curvature isotropic}.
\]

For $p\in M$ we define
\[
V_p=\gamma_\phi(p)={\mathbb R}\phi(p) \oplus \phi_{*} (T_pM).
\]
Then
\[
 M\times \R{{n+2}} =V\oplus V^{\perp}.
\]
The bundle $V^{\perp}$ is  flat by assumption, but $V$ is flat, too.
To see this,
we consider the second fundamental form $\beta$ of $V$:
\[
\beta: TM\times V \rightarrow V^{\perp}, \
\beta(X,s)=(D_Xs)^{V^{\perp}}.
\]
Then
\begin{equation}
<R^V(X,Y)s,t>=<\beta(X,s),\beta(Y,t)>-<\beta(Y,s),\beta(X,t)>.\label{gaussV}
\end{equation}
But for \(X,Y \in TM\)
\begin{eqnarray}
\beta(X,\phi)&=&0 \label{phikern}\\
\beta(X,\phi_{*}Y)&=&II_\phi(X,Y),\label{2ff}
\end{eqnarray}
where $II_\phi$ denotes the second fundamental form of
$\phi$. Since $M$ and $S$ have the same  curvature, this together with the
Gauss equation
implies, that the right-hand side  of (\ref{gaussV}) is zero.

Now  $V$ and $V^{\perp}$ are the pull-backs of the canonical bundle $\nu$
and its orthogonal complement $\nu^{\perp}$ over  $G_{m+1}(\R{{n+2}})$.
The tangent bundle of the Grassmannian can be
identified with Hom($\nu^{\perp},\nu$) and for tangent vector fields $X,Y,Z$
the curvature
tensor of $G_{m+1}(\R{{n+2}})$ is given by
\[
(R(X,Y)Z)(s)=R^{\nu}(X,Y)Z(s) - Z(R^{\nu^{\perp}}(X,Y)s),\; s\in
\Gamma(\nu^\perp).
\]
The flatness of $V$ and $V^\perp$ therefore implies that $\gamma_\phi$ is
curvature
isotropic.

If we canonically identify
\[
T_{\gamma_\phi}G_{m+1}(\R{{n+2}})={\text Hom}(V_p,V_p^{\perp}),
\]
then
\[
d_{p}\gamma_\phi(X)(s)=\beta(X,s).
\]
If we assume that the kernel of the second fundamental form of $\phi$ is
trivial or,
equivalently, that $f=\gamma_\phi$ is immersive, then we see from
(\ref{phikern}), (\ref{2ff}) that $\phi(p)$ is characterized uniquely up to
sign as a unit
vector in the kernel of $d_p\gamma_\phi(T_pM).$ This allows the reconstruction
of
$\phi$ from the associated curved flat.

This example can be worked out similarly for isometric immersions between
hyperbolic
spaces of equal curvature. The ambient space then is the Lorentzian ${\mathbb
R}^{n+2}_1$,
and the Gauss map is a space-like map into the $(m+1)$-planes of index 1.

Conversely, we can construct isometric immersions between space forms
of equal curvature from curved flats into the appropriate Grassmannians:
if
\[
f: M^m \to G_{m+1}(\R{{2m+1}})
\]
or
\[
f: M^m \to G_{m+1}(\Lor)
\]
is a space-like regular curved flat, then $A_1(T_{p}M)$ is $m$-dimensional for
all $p$.
Moreover, if it is a Cartan subalgebra (see Remark~\ref{rem:transitivity}
below),
then there exists a
gauge $H:M\to K$, where $K$ is either ${\bf O}(m+1)\times{\bf O}(m)$ or
${\bf O}(1,m)\times{\bf O}(m)$, so that the gauged frame $\tilde{F}=FH$ has
\[
\tilde{F}^{-1}\d\tilde{F}=
\left(
\begin{array}{c|c|cccc}
&\mp \sigma^T&&&&\\
\hline
&&-\beta_1&&&\\
\sigma&\omega&&\ddots&&\\
&&&&-\beta_m&\\
\hline
&*&&\eta&&
\end{array}
\right).
\]
This follows from Remark~\ref{rem:transitivity} below together with the fact
that the $\p$-part
of the above matrix, $\tilde A_{1}$, represents a normal form of space-like
Cartan subalgebras
in $\p$. Thus we see that  $\tilde{F}(p)e_0 \in \mbox{ker}\tilde{A}_1(T_pM)$.
Now consider
the map $\phi=\tilde{F}e_0$ which has its image in  $S^{2m}$ or $H^{2m}$. If
$\phi$ is
immersive then $\tilde{F}$ is an adapted frame for $\phi$ and (\ref{gaussV})
implies that $\phi$ is an isometric immersion between space forms of equal
curvature and flat
normal bundle.
\end{example}

\begin{example}\label{ex:isomnonequal} Let $M^m$ be of constant curvature $c\in
]0,1[$,
and let
\[
\phi_0: M^m \rightarrow S^n(1)
\]
be an isometric immersion with flat normal bundle. Note that the
latter condition is automatic if $n=2m-1$, see \cite{Moore}. We embed
$S^n(1)$ as an umbilic hypersurface into a sphere
$S=S^{n+1}(c) \subset \R{{n+2}}$ of curvature $c,$ and  obtain an isometric
immersion $\phi$ between constant curvature spaces of {\em equal}
curvatures as in the the preceeding example.

This construction is not limited to the above curvature relations. Consider an
$m$-dimensional space form M of curvature
$c$ isometrically immersed with flat normal bundle into an $n$-dimensional
space
form of curvature $\tilde c\neq c\neq 0$. The latter can always be realized as
an
umbilic hypersurface of a hyperquadric $S$ again of curvature $c$ in some
Euclidean or pseudo-Euclidean $(n+2)$-space. This produces examples of
curvature
isotropic maps into the Grassmanianns  given in the following list.

\renewcommand{\arraystretch}{2.0}
\begin{center}\begin{tabular}{|c|c|c|c|}
\hline
curvatures&S&ambient space&Grassmannian\\
\hline
$\tilde{c}>c>0$ & $S^{n+1}(c)$ & \R{{n+2}}&%
$\frac{{\bf O}(n+2)}{{\bf O}(m+1)\times {\bf O}(n+1-m)}$\\
\hline
$\tilde{c}<c>0$ & $S^{n+1}_1(c)$ & $ {\mathbb R}_1^{n+2}$ &%
$\frac{{\bf O}(1,n+1)}{{\bf O}(1+m)\times {\bf O}(1,n-m)}$\\
\hline
$\tilde{c}>c<0$ & $H^{n+1}(c)$ &${\mathbb R}_1^{n+2}$&%
$\frac{{\bf O}(1,n+1)}{{\bf O}(1,m)\times {\bf O}(n+1-m)}$\\
\hline
$\tilde{c}<c<0$ & $ H^{n+1}_1(c)$ & ${\mathbb R}_2^{n+2}$&%
$\frac{{\bf O}(2,n)}{{\bf O}(1,m)\times {\bf O}(1,n-m)}$\\
\hline
\end{tabular}\end{center}

The isometric immersions between space forms of equal curvature that arise in
the above manner carry a distinguished umbilic normal field, namely that of
the ambient umbilic hypersurface of $S$.
Conversely, consider an isometric immersion $\phi$ with flat normal bundle and
immersive Gauss map between space forms of equal curvature. Assume that we are
in
the critical codimension $n=2m-1$. Then we can show the existence of a
distinguished umbilic normal field $\xi$, and prove that the constancy of
$|\xi|$
is necessary and sufficient for $\phi$ to factor through an umbilic
hypersurface
of $S$.

A different approach for the construction of isometric immersions of space
forms
of distinct curvatures is presented in \cite{FP}.
\end{example}

\begin{example}\label{ex:isotherm}
 A curved 2-flat in the Grassmannian $G^+_3(\Lorentz{5})$ of space-like
3-planes in
Lorentzian 5-space has a lift
\[
F: M^2 \rightarrow {\bf O}(1,4).
\]
The Lorentz group acts on the sphere $S^4_1:=\{v\in\Lorentz{5}\;|<v,v>=1\}$,
and an
appropriate orbit of $F$ gives a special Ribeaucour sphere congruence in $S^3$
enveloping
a  pair of isothermic surfaces. See \cite{BJPP}, p.4.
\end{example}

\section{Geometry of curved flats}\label{sec:geometry}

We have seen that a curved flat $f:M\to N$ is at each point tangential to {\em
flats} in the
symmetric space $N$. As we shall see below this implies (under certain
non-degeneracy
assumptions) that a curved flat is intrinsically flat (where  the metric on $N$
is always the one induced by the Killing form of $\mathfrak g$).
In this sense curved flats
can be regarded as analogues of developable surfaces in Euclidean space.
Moreover, if the curved flat is the Gauss map of an isometric immersion then
the developing isometry induces principal curvature coordinates on $M$. In this
situation
the flat metric has been used to show certain non-existence results
for isometric immersions \cite{P}.

\begin{mytheorem}\label{th:flat}
Let $f:M \rightarrow N$ be a regular curved flat whose tangent spaces
are all conjugate under the action of the connected component of the isotropy
group of $N$ .
Then $M$ with the induced metric
\(f^{*}g_{N}\) is (intrinsically) flat.
\end{mytheorem}

\begin{remark}\label{rem:transitivity}
The transitivity assumption of the Theorem will hold provided that all tangent
spaces of the
curved flat are Cartan subalgebras: this follows from the fact that the
connected components
of the space of all Cartan subalgebras in $\p$ are precisely the orbits under
(the connected
component of) the isotropy group $K_0$. Recall that a subspace $\a$ of $\p$
is a {\em Cartan subalgebra} if  $\a$ is maximal abelian and consists of
semi-simple elements. All Cartan subalgebras have the same dimension
\cite{Berger},
say $r$,
which defines the rank of the (pseudo-Riemannian) symmetric space $G/K$. Now
consider
\[
 \Sigma = \{\a\subset\p^{\C}\,;\, \a\; \text{abelian and}\, \dim_{\C}\a =
r\,\}\subset G_r(\p^{\C})
\]
which is a subvariety (defined over $\Rr$) of the Grassmanian of $r$-planes in
$\p^{\C}$.
The set of all Cartan subalgebras
\[
\widetilde{\Sigma}=\{\a\in\Sigma\,;\,\a\; \text{Cartan}\,\}
\]
is open in $\Sigma$ since it is easy to show that a maximal abelian $\a$
is Cartan iff the  Killing form is non-degenerate
on $\a$ and the set of maximal abelian $\a\in\Sigma$ is open in $\Sigma$.
Moreover, $K^{\C}$ acts transitively on $\widetilde{\Sigma}$ with
stabilizer $N(\a_{0})=\{ k\in K^{\C}\,;\, \Ad K^{\C}(\a_{0})=\a_{0}\,\}$,
$\a_{0}\in\widetilde{\Sigma}$,
so that all points in $\widetilde{\Sigma}\subset\Sigma$ are smooth and
$\widetilde{\Sigma}=
K^{\C}/N(\a_{0})$. Let
\[
\Sigma(\Rr) = \{\a\subset\p\,;\,\a\;\text{abelian and}\,\dim_{\Rr}\a =
r\,\}\subset G_r(\p)
\]
be the $\Rr$-points of the variety $\Sigma$. Then the set of real Cartan
subalgebras in $\p$
is
\[
\widetilde{\Sigma}_{\Rr} = \widetilde{\Sigma}\cap\Sigma(\Rr)
\]
and as above $\widetilde{\Sigma}_{\Rr}$ is open in $\Sigma(\Rr)$. Moreover,
\cite{Whitney},
$\widetilde{\Sigma}_{\Rr}$ is a smooth real manifold of
$\dim_{\Rr}\widetilde{\Sigma}_{\Rr}=
\dim_{\C}\widetilde{\Sigma}$ with only finitely many connected components.
To see that each $K_0$-orbit $\mathcal{O}$
is a connected component of $\widetilde{\Sigma}_{\Rr}$ it suffices to show that
 $\mathcal{O}$ is open in  $\widetilde{\Sigma}_{\Rr}$.
Since $\mathcal{O}=K_0/N(\a_{0})$ where $N(\a_0)$ is
the normalizer
of some $\a_0\in\widetilde{\Sigma}_{\Rr}$ we obtain
\[
\dim_{\Rr}\mathcal{O}=\dim_{\C}\widetilde{\Sigma}=\dim_{\Rr}\widetilde{\Sigma}_{\Rr}
\]
which shows that $\mathcal{O}$ is open in $\widetilde{\Sigma}_{\Rr}$.

In the  above discussion we saw that being a Cartan subalgebra is an open
condition inside
rank dimensional abelian subspaces. When we integrate the curved
flat equations in Section~\ref{sec:integrate} we can always assure this locally
by choosing the
right initial conditions for the Lax flows.

Of course, for Riemannian symmetric spaces N it suffices to have $\dim
M\,=\,\rank N$ for the
tangent spaces of a regular curved flat $f:M\to N$ to be Cartan.
\end{remark}

\begin{proof}[Proof of the Theorem]
We argue locally at \(p_{0} \in M\), and may hence assume the
existence of a framing \(F:M \rightarrow G\) of \(f\), such that
as in (\ref{decomp})
\begin{eqnarray}
F^{-1}\d F&=&A=A_0+A_1 \,,
\\
F(p_{0})&=&e \in G\,.
\end{eqnarray}
Then
\[
\mathfrak{a}_p:=f_{{*}}(T_{p}M)=A_{1}(T_{p}M) \subset \mathfrak{p}
\]
is abelian by Lemma~\ref{lem:abel}.
We change \(F\) by a gauge transformation using the transitivity assumption.
Let \(H: M \rightarrow K\) be such that
\(\Ad_{H(p)}(\mathfrak{a}_{p})=\mathfrak{a}_{p_0}\), and put \(\tilde{F}:=FH\).
Then taking the $\mathfrak{p}$-part of the Maurer-Cartan equation
(\ref{MC}) for the corresponding
\[
\tilde{A}_{1}=\Ad(H)A_{1}:TM \rightarrow \mathfrak{a}_{p_0}
\]
we find
\begin{equation}
\d\tilde{A}_{1}+ 2 [\tilde{A}_{0}\wedge\tilde{A}_{1}]=0.\label{eq:tildeMC}
\end{equation}
Clearly \(\d\tilde{A}_{1}\) has values in \(\mathfrak{a}_{p_0}\),
while
\[
<[\mathfrak{k},\mathfrak{a}_{p_0}],\mathfrak{a}_{p_0}>=<\mathfrak{k},%
[\mathfrak{a}_{p_0},\mathfrak{a}_{p_0}]>=0
\]
implies that \([\tilde{A}_{0}\wedge\tilde{A}_{1}]\) has values in
\(\mathfrak{a}_{p_0}^{\perp}\). Since $f$ was assumed to be regular
$\mathfrak{a}_{p_0}$ and $\mathfrak{a}_{p_0}^{\perp}$ are complementary
subspaces.
\\
Therefore
\begin{equation}
\d\tilde{A}_{1}=0, \label{eq:closed}
\end{equation}
and so \(\tilde{A}_1=\d\psi\) for some local diffeomorphism
\(\psi:M\rightarrow \mathfrak{a}_{p_0}\).
Finally
\begin{eqnarray*}
<f_{{*}}X,f_{{*}}Y>&=&<\pi_{{*}}\tilde{F}_{{*}}X,\pi_{{*}}\tilde{F}_{{*}}Y>
\\
&=&<\tilde{A}_{1}(X),\tilde{A}_{1}(Y)>_{\mathfrak{a}_{p_0}}
 \\                      &=&<d\psi(X),d\psi(Y)>_{\mathfrak{a}_{p_0}}
\end{eqnarray*}
shows that the metric is flat.
\end{proof}

\begin{remark}
Let $f:M\to N$ be a regular curved flat such that all its tangent spaces are
Cartan
subalgebras. As in the proof let $\a=f_{{*}}(T_{p_{0}}M)$ be a fixed Cartan
subspace
in $\p$ and $\psi:M\to\a$ the developing isometry. If $\alpha_1,\dots,\alpha_m$
are
a set of simple roots on $\a$ then $\alpha_{i}\circ\psi$ give canonical
coordinates
on $M$. If the curved flat arises from an isometric immersion between space
forms
like at the end of Example \ref{ex:isomequal} then one can easily check that
these coordinates are the principal curvature coordinates for the immersion.
This provides a rather nice geometric interpretation of the developing isometry
of a curved flat.
\end{remark}

\section{The spectral parameter}\label{sec:spectral}

In this section we reformulate the curved flat equations \eqref{MC},
\eqref{flat} as an integrability condition involving an additional -
{\em spectral} - parameter. If the curved flat happens to be the Gauss
map of an isometric immersion (Example~\ref{ex:isomequal}) this parameter
scales the second fundamental form of the immersion and thus gives
a non-trivial deformation. The existence of such a non-trivial deformation is
indicative that curved flats are solutions to a certain integrable
system. We will discuss this aspect in more detail in
Section~\ref{sec:integrate}.

Let $A_0$ and $A_1$ be $\mathfrak k$ and $\mathfrak p$-valued 1-forms and
consider
the 1-parameter family ({\em loop}) of $\mathfrak g$-valued 1-forms
\begin{equation}
A^\mu := A_0+ \mu A_1  \label{loop}
\end{equation}
for real $\mu\neq 0$. Then we have
\begin{mylemma}\label{lem:loop}
The following are  equivalent statements:
\begin{enumerate}
\item
$A_0$ and $A_1$ satisfy the curved flat equations (\ref{MC}) and (\ref{flat});
\item
$A_0$ and $\mu A_1$ satisfy the curved flat equations (\ref{MC}) and
(\ref{flat}) for
some $\mu\neq 0$;
\item
$A^\mu$ solves the Maurer-Cartan equation
\[
\d A^{\mu} + \frac{1}{2}[A^{\mu}\wedge A^{\mu}]=0
\]
identical in $\mu$.
\end{enumerate}
\end{mylemma}
The proof follows immediately by comparing coefficients at powers of $\mu$.
Thus we see that curved flats come in 1-parameter families and that their
integrability equations can be rewritten
as a zero-curvature condition involving and auxiliary parameter.

According to \cite{FP}, an isometric immersion $f:M(c)\to\tilde{M}(\tilde{c})$
with flat normal bundle
between space forms
of distinct curvatures also comes in a 1-parameter family --- let us call this
parameter $\lambda$. As explained in
Example~\ref{ex:isomnonequal}, their Gauss
map is a curved flat, which by the above also embeds in a 1-parameter family.
It can be shown that the two families correspond where the relation between the
two parameters is given by
\[
\mu=\frac{-\sqrt{c}}{2\sqrt{1-c}}(\lambda-\lambda^{-1})\,.
\]

\section{Integration of the Curved Flat Equations}\label{sec:integrate}
Using the reformulation of the curved flat equations (Lemma~\ref{lem:loop}) as
a loop of flat $1$-forms we will now apply the finite type integration scheme
to construct solutions in terms of a hierarchy of finite dimensional
commuting ODE in Lax form.
Let $N = G/K$ be a semisimple symmetric space,
\[
\mathfrak{g}=\mathfrak{k}+\mathfrak{p}
\]
the symmetric decomposition of $G/K$, and denote by
$\sigma:\mathfrak{g}\rightarrow \mathfrak{g}$ the corresponding involution.
We put
\begin{eqnarray}
\Lambda\mathfrak{g_\sigma}&:=&\{\xi:{\mathbb R}^*\to \mathfrak{g}
\mbox{ Laurent polynomial in }\mu \;|\;\sigma\xi(\mu)=\xi(-\mu)\},
\end{eqnarray}
and split this as
\[
\Lambda\mathfrak{g_\sigma}=\Lambda^+\mathfrak{g_\sigma}
\oplus\Lambda^-_*\mathfrak{g_\sigma}.
\]
Here $\Lambda^+\mathfrak{g_\sigma}$ denotes the polynomial loops in $\mu$,
while those in
$\Lambda^-_*\mathfrak{g_\sigma}$ are polynomials in $\mu^{-1}$ that vanish at
infinity. Let
$\pi_\pm$ denote the corresponding projections. Then
\[
R:= \frac{1}{2}(\pi_+ - \pi_-):\Lambda\mathfrak{g_\sigma}
\to \Lambda\mathfrak{g_\sigma}
\]
defines an R-matrix, and we can apply the results of \cite{FPPS}, Section 2,
and
\cite{BFPP}, Section 3.

Let $d>0$ be odd, and let $V$ be an Ad$G$-equivariant polynomial vector field
on
$\mathfrak{g}$. We decompose this into its homogeneous parts of degrees
$r_1,\ldots,r_m$,
\[
V=V^{(1)}+\ldots + V^{(m)},
\]
and define for $\xi\in\Lambda\mathfrak{g}_\sigma$
\[
\tilde{V}(\xi)(\mu):=\sum_{i=1}^m (-1)^{r_i+1}\mu^{1-dr_i}V^{(i)}(\xi(\mu)).
\]

If $V$ is $(-\sigma)$-equivariant
\[
-\sigma V(\xi)=V(-\sigma \xi),
\]
then $\tilde{V}$ is ad-equivariant on $\Lambda\mathfrak{g}_\sigma$, and
induces a vector field

\[
X_V(\xi)=[ \xi,(R+\frac{1}{2})\tilde{V}(\xi)],
\]

on $\Lambda\mathfrak{g_\sigma}$, which is in fact tangential to the finite
dimensional subspace
\[
\Lambda_d := \{\xi\in\Lambda^+\mathfrak{g_\sigma}|
\xi=\sum_{k=0}^{d}\mu^k\xi_k\}.
\]
Any $k$ such vector fields $C^\infty$-commute. Natural examples are given by
the gradients ${\rm grad}\, v=V$  of  Ad$G$- and $(-\sigma)$-invariant
polynomials $v$ (with respect to an Ad$G$-invariant inner product on
$\mathfrak{g}$).
If N is of rank $k$, then the existence of $k$ independent such
functions is guaranteed, provided N has the surjection property, see
\cite{BFPP} Section 8,
which is the case for all classical symmetric spaces and group manifolds.

We now are in the position to state the basic integration theorem for curved
flats:
\begin{mytheorem}
Given $\Ad G$- and $(-\sigma)$-equivariant polynomial vector fields
$V_1,\ldots,V_k$ on $\mathfrak{g}$, and
$\xi_0 \in \Lambda_d$, for some odd $d>0$, the system
\[
\frac{\partial \xi}{\partial x_j}=[\xi,(R+\frac{1}{2})\tilde{V_j}(\xi)]
\]
has a unique local solution $\xi: \R{k} \supset U \rightarrow \Lambda_d$
with $\xi(0)=\xi_0$. For this
\[
A^\mu = \sum_j(R+\frac{1}{2})\tilde{V_j}(\xi)dx^j
\]
is of the form (\ref{loop}). It satisfies (\ref{MC}), (\ref{flat}) for all
$\mu$,
and therefore induces a curvature isotropic map $f: \R{k} \supset U \rightarrow
N$.
\end{mytheorem}
\begin{remark}
If the rank of $G/K$ is $k$ then we can choose the initial condition $\xi_0$
in such a way that $(R+\frac{1}{2})\tilde{V_j}(\xi_0)$, $j=1,\dots ,k$ span
a Cartan subalgebra in $\p$. But this is an open condition inside
$k$-dimensional
abelian subspaces and thus
each tangent space of the resulting curved flat will be a Cartan subalgebra
for $p\in U$ (c.f. Remark~\ref{rem:transitivity}).
\end{remark}
\begin{remark}
If $G/K$ is a Grassmannian --- which is the case for the Gauss maps of
isometric immersions ---
the vector fields $V_i$ are simply given by even powers of off-blockdiagonal
matrices
$\xi^{2i}$.
\end{remark}


\begin{thebibliography}{99}

\bibitem{Berger}
 M. Berger, {\em  Les espaces sym\'etriques non compacts},
 Ann. Ec. Norm. Sup. {\bf 74} (1957), 85-177.

\bibitem{BJPP}
F. Burstall, U. Hertrich-Jeromin, F. Pedit, U. Pinkall,
{\em Curved flats and isothermic surfaces}, SFB 288 Preprint No. 132 (1994)

\bibitem{BFPP}
F.E. Burstall, D. Ferus, F. Pedit, U. Pinkall,
{\it Harmonic tori in symmetric spaces and commuting %
Hamiltonian systems on loop algebras},
Ann. of Math. {\bf 138} (1993),173-212.

\bibitem{FP}
D. Ferus, F. Pedit,
{\it Isometric immersions of space forms and soliton theory},
SFB 288 Preprint No. 154 (1995).

\bibitem{FPPS}
D. Ferus, F. Pedit, U. Pinkall, I. Sterling,
{\it Minimal tori in $S^4$},
J. Reine Angew. Math. {\bf 429} (1992), 1-47.


\bibitem{Moore}
J.D. Moore,
{\it Isometric immersions of space forms into space forms},
Pacific J. Math. {\bf 40} (1972),157-166.


\bibitem{P}
F. Pedit,
{\it A non-immersion theorem for space forms},
Comment. Math. Helvetici {\bf 63} (1988), 672-674.


\bibitem{UJ}
U. Hertrich-Jeromin,
{\it Conformally flat hypersurfaces and integrable systems},
Thesis, TU-Berlin (1994).


\bibitem{UJ1}
U. Hertrich-Jeromin,
{\it  On  conformally flat hypersurfaces, curved flats and cyclic systems},
in preparation.

\bibitem{Whitney}
H. Whitney,
{\it Elementary structure of real algebraic varieties}, Ann. of Math. {\bf 66}
(1957), 545-556.

\end{thebibliography}
\end{document}